\documentclass[12pt, reqno]{amsart}
\usepackage{amsmath, amsthm, amscd, amsfonts, amssymb, graphicx}
\usepackage[bookmarksnumbered,plainpages]{hyperref}

\newtheorem{theorem}{Theorem}[section]

\theoremstyle{definition}

\theoremstyle{remark}

\numberwithin{equation}{section}

\begin{document}
\title[Some inequalities for $(\alpha, \beta)$-normal operators]{Some inequalities for $(\alpha, \beta)$-normal operators in Hilbert
spaces}
\author[S.S. Dragomir, M.S. Moslehian]{S. S. Dragomir$^{1}$ and M. S. Moslehian$^{2}$}
\address{$^{1}$ School of Computer Science and Mathematics, Victoria
University, P. O. Box 14428, Melbourne City, Victoria 8001, Australia.}
\email{sever.dragomir@vu.edu.au}
\address{$^{2}$ Department of Mathematics, Ferdowsi University of Mashhad,
P.O. Box 1159, Mashhad 91775, Iran\\
Centre of Excellence in Analysis on Algebraic Structures (CEAAS), Ferdowsi
University of Mashhad, Iran.}
\email{moslehian@ferdowsi.um.ac.ir and moslehian@ams.org}

\begin{abstract}
An operator $T$ acting on a Hilbert space is called $(\alpha ,\beta )$%
-normal ($0\leq \alpha \leq 1\leq \beta $) if
\begin{equation*}
\alpha ^{2}T^{\ast }T\leq TT^{\ast }\leq \beta ^{2}T^{\ast }T.
\end{equation*}%
In this paper we establish various inequalities between the operator norm
and its numerical radius of $(\alpha ,\beta )$-normal operators in Hilbert
spaces. For this purpose, we employ some classical inequalities for vectors
in inner product spaces.
\end{abstract}

\maketitle

\section{Introduction}

An operator $T$ acting on a Hilbert space $({\mathcal{H}};\langle \cdot
,\cdot \rangle )$ is called $(\alpha ,\beta )$\textit{-normal} ($0\leq
\alpha \leq 1 \leq \beta $) if
\begin{equation*}
\alpha ^{2}T^{\ast }T\leq TT^{\ast }\leq \beta ^{2}T^{\ast }T.
\end{equation*}%
Then
\begin{equation*}
\alpha ^{2}\langle T^{\ast }Tx,x\rangle \leq \langle TT^{\ast }x,x\rangle
\leq \beta ^{2}\langle T^{\ast }Tx,x\rangle \,,
\end{equation*}%
whence
\begin{equation}
\alpha \Vert Tx\Vert \leq \Vert T^{\ast }x\Vert \leq \beta \Vert Tx\Vert ,
\label{nor}
\end{equation}%
for all $x\in {\mathcal{H}}$, namely, both $T$ majorizes $T^*$ and $T^*$
majorizes $T$. A seminal result of R.G. Douglas \cite{DOU} (majorization
lemma) says that an operator $T \in B({\mathcal{H}})$ majorizes an operator $%
S\in B({\mathcal{H}})$ if any one of the following equivalent statements
holds:

(i) the range space $\mathrm{ran}(T)$ of $T$ is a subset of $\mathrm{ran}(S)$%
;

(ii) $TT^{\ast }\leq \lambda ^{2}SS^{\ast }$;

(iii) there exists an operator $R\in B({\mathcal{H}})$ such that $T=SR$.%
\newline
Furthermore, $R$ is the unique operator satisfying

(a) $\Vert R\Vert =\inf \{\lambda :TT^{\ast }\leq \lambda SS^{\ast }\}$;

(b) $\ker (T)=\ker (R)$;

(c) $\mathrm{ran}(R)$ is a subset of the closure $\mathrm{ran}(S^*)^-$ of $%
\mathrm{ran}(S^*)$.

\noindent Analogues of Douglas' majorization lemma for Banach space
operators were studied by M.R. Embry \cite{EMB} (see also \cite{BOU1}). A
discussion of the duality between the properties of majorization and range
inclusion can be found in \cite{BAR}.

\noindent Using the result of Douglas, we observe that $T$ is $(\alpha
,\beta )$-normal if and only if $\mathrm{ran}(T)=\mathrm{ran}(T^{\ast })$,
or, equivalently, $\ker (T)=\ker (T^{\ast })$. It is therefore obvious that
invertible, normal and hyponormal operators are $(\alpha ,\beta )$-normal
for some appropriate values of $\alpha $ and $\beta $. The matrix
{\scriptsize $\left[
\begin{array}{cc}
1 & 0 \\
1 & 1%
\end{array}%
\right] $} in $B({\mathbb{C}}^{2})$ is an $(\alpha ,\beta )$-normal with $%
\alpha =\sqrt{(3-\sqrt{5})/2}$ and $\beta =\sqrt{(3+\sqrt{5})/2}$, which is
neither normal nor hyponormal. There are some results which can be applied
to our notion in the literature. For instance, one can deduce from \cite[%
Lemma 1]{BOU2} that if $z$ is an eigenvalue of $T$ and $z$ belongs to the
topological boundary of the numerical range of $T$, then $T-z$ is $(\alpha
,\beta )$-normal for some $\alpha $ and $\beta $. There are also some
interesting questions in linear algebra concerning $(\alpha ,\beta )$%
-normality, see \cite{MOS}.

\noindent Another characterization is that $T$ is $(\alpha ,\beta )$-normal (%
$0 < \alpha \leq 1 \leq \beta $) if and only if there are operators $S_1,
S_2 \in B({\mathcal{H}})$ such that $T=T^*S_1$ and $T=S_2T^*$. Moreover, $%
S_1, S_2$ can be chosen in such a way that
\begin{equation*}
\|S_1\|=\inf\{\beta \geq 1: TT^* \leq \beta T^*T\}, \qquad
\|S_2\|=\sup\{\alpha>0: \alpha T^*T \leq TT^*\}\,.
\end{equation*}
Let $T$ be an $(\alpha ,\beta )$-normal operator on a (not necessarily
finite dimensional) Hilbert space ${\mathcal{H}}$. Using the fact that $%
ker(T^*)^\perp=\mathrm{ran}(T)^-$, we observe that ${\mathcal{H}}=\ker(T)
\oplus \mathrm{ran}(T)^-$. Hence $T$ can be represented as a block matrix
{\scriptsize $\left[
\begin{array}{cc}
0 & 0 \\
0 & C%
\end{array}%
\right]$}, where $C: \mathrm{ran}(T)^- \to \mathrm{ran}(T)^-$ has zero
kernel. We can define the pseudo-inverse of $T$, denoted by $T^+$, to be the
operator on ${\mathcal{H}}$, which is zero on $\mathrm{ran}(T)^\perp$, and
is the inverse to $C$ on $\mathrm{ran}(T)^-$. It is easy to see that $T^+$
is closed if and only if $\mathrm{ran}(T)$ is closed. The operator
pseudo-inverse is a powerful tool in applied mathematics; cf. \cite{B-R}.

Let $({\mathcal{H}};\langle \cdot ,\cdot \rangle )$ be a complex Hilbert
space. The \textit{numerical radius} $w(T)$ of an operator $T$ on ${\mathcal{%
H}}$ is given by
\begin{equation}
w(T)=\sup \{|\langle Tx,x\rangle |,\Vert x\Vert =1\}.  \label{1.1}
\end{equation}%
Obviously, by (\ref{1.1}), for any $x\in \mathcal{H}$ one has
\begin{equation}
|\langle Tx,x\rangle |\leq w(T)\Vert x\Vert ^{2}.  \label{1.2}
\end{equation}

\noindent It is well known that $w(\cdot)$ is a norm on the Banach algebra $%
B({\mathcal{H}})$ of all bounded linear operators. Moreover, we have
\begin{equation*}
w(T)\leq \Vert T\Vert \leq 2w(T)\qquad (T\in B(\mathcal{H})).
\end{equation*}%
For other results and historical comments on the numerical radius see \cite%
{GR}.

In this paper, we establish various inequalities between the operator norm
and its numerical radius of $(\alpha ,\beta )$-normal operators in Hilbert
spaces. For this purpose, we employ some classical inequalities for vectors
in inner product spaces due to Buzano, Dunkl--Williams, Dragomir--S\'{a}%
ndor, Goldstein--Ryff--Clarke and Dragomir.

\section{Inequalities Involving Numerical Radius}

In this section we study some inequalities concerning the numerical radius
and norm of $(\alpha ,\beta )$-normal operators. Our first result reads as
follows, see also \cite{SSD5}:

\begin{theorem}
\label{t.2b} Let $T\in B({\mathcal{H}})$ be an $(\alpha ,\beta )$-normal
operator. Then
\begin{equation}
(\alpha ^{2r}+\beta ^{2r})\Vert T\Vert ^{2}\leq \left\{
\begin{array}{ll}
2\beta ^{r}w(T^{2})+r^{2}\beta ^{2r-2}\Vert \beta T-T^{\ast }\Vert ^{2}, &
\text{ if \ }r\geq 1, \\
&  \\
2\beta ^{r}w(T^{2})+\Vert \beta T-T^{\ast }\Vert ^{2}, & \text{if \ }r<1.%
\end{array}%
\right.  \label{2.5b}
\end{equation}
\end{theorem}

\begin{proof}
We use the following inequality for vectors in inner product spaces due to
Goldstein, Ryff and Clarke \cite{GRC}:
\begin{equation}
\Vert a\Vert ^{2r}+\Vert b\Vert ^{2r}-2\Vert a\Vert ^{r}\Vert b\Vert
^{r}\cdot \frac{\text{Re}\langle a,b\rangle }{\Vert a\Vert \,\Vert b\Vert }%
\leq \left\{
\begin{array}{ll}
r^{2}\Vert a\Vert ^{2r-2}\Vert a-b\Vert ^{2} & \text{ if \ }r\geq 1, \\
&  \\
\Vert b\Vert ^{2r-2}\Vert a-b\Vert ^{2} & \text{if \ }r<1,%
\end{array}%
\right.  \label{2.6b}
\end{equation}%
provided $r\in \mathbb{R}$ and $a,b\in H$ with $\Vert a\Vert \geq \Vert
b\Vert .$

Suppose that $r\geq 1$. Let $x\in H$ with $\Vert x\Vert =1$. Noting to (\ref%
{nor}) and applying (\ref{2.6b}) for the choices $a=\beta Tx$, $b=T^{\ast }x$
we get
\begin{multline}
\Vert \beta Tx\Vert ^{2r}+\Vert T^{\ast }x\Vert ^{2r}-2\Vert \beta Tx\Vert
^{r-1}\Vert \,\Vert T^{\ast }x\Vert ^{r-1}\,\text{Re}\langle \beta
Tx,T^{\ast }x\rangle  \label{2.7b} \\
\leq r^{2}\Vert \beta Tx\Vert ^{2r-2}\Vert \beta Tx-T^{\ast }x\Vert ^{2}
\end{multline}%
for any $x\in H,$ $\Vert x\Vert =1$ and $r\geq 1.$ Using (\ref{nor}) and (%
\ref{2.7b}) we get
\begin{multline}
(\alpha ^{2r}+\beta ^{2r})\Vert Tx\Vert ^{2r}  \label{2.7bb} \\
\leq 2\beta ^{r}\Vert Tx\Vert ^{r-1}\Vert T^{\ast }x\Vert ^{r-1}|\langle
T^{2}x,x\rangle |+r^{2}\beta ^{2r-2}\Vert Tx\Vert ^{2r-2}\Vert \beta
Tx-T^{\ast }x\Vert ^{2}.
\end{multline}%
Taking the supremum in (\ref{2.7bb}) over $x\in H,$ $\Vert x\Vert =1,$ we
deduce
\begin{equation*}
(\alpha ^{2r}+\beta ^{2r})\Vert T\Vert ^{2r}\leq 2\beta ^{r}\Vert T\Vert
^{2r-2}\Vert T^{\ast }\Vert ^{r-1}w(T^{2})+r^{2}\beta ^{2r-2}\Vert T\Vert
^{2r-2}\Vert \beta T-T^{\ast }\Vert ^{2},
\end{equation*}%
which is the first inequality in (\ref{2.5b}). If $r<1$, then one can
similarly prove the second inequality in (\ref{2.5b}).
\end{proof}

\begin{theorem}
Let $T\in B({\mathcal{H}})$ be an $(\alpha ,\beta )$-normal operator. Then
\begin{equation}
w(T)^{2}\leq \frac{1}{2}\left[ \beta \Vert T\Vert ^{2}+w(T^{2})\right] .
\label{dw}
\end{equation}
\end{theorem}

\begin{proof}
The following inequality is known in the literature as the \textit{Buzano
inequality} \cite{BUZ}:
\begin{equation}
|\langle a,e\rangle \langle e,b\rangle |\leq \frac{1}{2}(\Vert a\Vert
\,\Vert b\Vert +|\langle a,b\rangle |),  \label{bu}
\end{equation}%
for any $a,b,e$ in ${\mathcal{H}}$ with $\left\Vert e\right\Vert =1.$

Let $x\in H$ with $\Vert x||=1$. Put $e=x,a=Tx,b=T^{\ast }x$ in (\ref{bu})
to get
\begin{multline*}
|\langle Tx,x\rangle \langle x,T^{\ast }x\rangle |\leq \frac{1}{2}(\Vert
Tx\Vert \,\Vert T^{\ast }x\Vert +|\langle Tx,T^{\ast }x\rangle |) \\
\leq \frac{1}{2}(\beta \Vert Tx\Vert ^{2}+|\langle T^{2}x,x\rangle |).
\end{multline*}%
Taking the supremum over $x\in H,$ $\Vert x\Vert =1,$ we obtain (\ref{dw}).
\end{proof}

\begin{theorem}
Let $T\in B({\mathcal{H}})$ be an $(\alpha ,\beta )$-normal operator and $%
\lambda \in {\mathbb{C}}$. Then
\begin{equation}
\alpha \Vert T\Vert ^{2}\leq w(T^{2})+\frac{2\beta \Vert T-\lambda T^{\ast
}\Vert ^{2}}{(1+|\lambda |\alpha )^{2}}.  \label{tt}
\end{equation}
\end{theorem}

\begin{proof}
Using the \textit{Dunkl--Williams inequality} \cite{D-W}
\begin{equation*}
\frac{1}{2}(\Vert a\Vert +\Vert b\Vert )\left\Vert \frac{a}{\Vert a\Vert }-%
\frac{b}{\Vert b\Vert }\right\Vert \leq \Vert a-b\Vert \qquad (a,b\in
H\setminus \left\{ 0\right\} )
\end{equation*}%
we get
\begin{equation*}
2-2\cdot \frac{\text{Re}\langle a,b\rangle }{\Vert a\Vert \Vert b\Vert }%
=\left\Vert \frac{a}{\Vert a\Vert }-\frac{b}{\Vert b\Vert }\right\Vert
^{2}\leq \frac{4\Vert a-b\Vert ^{2}}{(\Vert a||+\Vert b\Vert )^{2}}\qquad
(a,b\in H\setminus \left\{ 0\right\} )
\end{equation*}%
whence
\begin{equation*}
\Vert a\Vert \Vert b\Vert \leq \frac{2\Vert a\Vert \,\Vert b\Vert \,\Vert
a-b\Vert ^{2}}{(\Vert a||+\Vert b\Vert )^{2}}+|\langle a,b\rangle |\qquad
(a,b\in H\setminus \left\{ 0\right\} ).
\end{equation*}%
Put $a=Tx$ and $b=\lambda T^{\ast }$ to get
\begin{equation*}
\Vert Tx\Vert \,\Vert T^{\ast }x\Vert \leq |\langle T^{2}x,x\rangle |+\frac{%
2\Vert Tx\Vert \,\Vert T^{\ast }x\Vert \,\Vert Tx-\lambda T^{\ast }x\Vert
^{2}}{(\Vert Tx\Vert +|\lambda |\,\Vert T^{\ast }x\Vert )^{2}}
\end{equation*}%
so that
\begin{multline}
\alpha \Vert Tx\Vert ^{2}\leq |\langle T^{2}x,x\rangle |+\frac{2\beta \Vert
Tx\Vert ^{2}\Vert Tx-\lambda T^{\ast }x\Vert ^{2}}{(\Vert Tx\Vert +|\lambda
|\alpha \,\Vert Tx\Vert )^{2}}  \label{tt*} \\
\leq |\langle T^{2}x,x\rangle |+\frac{2\beta \Vert (T-\lambda T^{\ast
})x\Vert ^{2}}{(1+|\lambda |\alpha )^{2}}.
\end{multline}%
Taking the supremum in (\ref{tt*}) over $x\in H,$ $\Vert x\Vert =1,$ we get
the desired result (\ref{tt}).
\end{proof}

\begin{theorem}
Let $T\in B({\mathcal{H}})$ be an $(\alpha ,\beta )$-normal operator and $%
\lambda \in {\mathbb{C}}\backslash \{0\}$. Then
\begin{equation}
\left[ \alpha ^{2}-\left( \frac{1}{|\lambda |}+\beta \right) ^{2}\right]
\Vert T\Vert ^{4}\leq w(T^{2}).  \label{drag0}
\end{equation}
\end{theorem}

\begin{proof}
We apply the following reverse of the quadratic Schwarz inequality obtained
by Dragomir in \cite{DRA}
\begin{equation}
(0\leq )\Vert a\Vert ^{2}\Vert b\Vert ^{2}-|\langle a,b\rangle |^{2}\leq
\frac{1}{\left\vert \lambda \right\vert ^{2}}\Vert a\Vert ^{2}\Vert
a-\lambda b\Vert ^{2}  \label{drag}
\end{equation}%
provided $a,b\in H$ and $\lambda \in {\mathbb{C}}\backslash \{0\}$.

Set $a=Tx,b=T^{\ast }x$ in (\ref{drag}), to get
\begin{multline*}
\alpha ^{2}\Vert Tx\Vert ^{4}\leq |\langle Tx,T^{\ast }x\rangle |^{2}+\frac{1%
}{|\lambda |^{2}}\Vert Tx\Vert ^{2}\Vert Tx-\lambda T^{\ast }x\Vert ^{2} \\
\leq |\langle T^{2}x,x\rangle |^{2}+\frac{1}{|\lambda |^{2}}\Vert Tx\Vert
^{2}(1+|\lambda |\beta )^{2}\Vert Tx\Vert ^{2}
\end{multline*}%
whence
\begin{equation}
\left[ \alpha ^{2}-\left( \frac{1}{|\lambda |}+\beta \right) ^{2}\right]
\Vert Tx\Vert ^{4}\leq |\langle T^{2}x,x\rangle |^{2}.  \label{drag1}
\end{equation}%
Taking the supremum in (\ref{drag1}) over $x\in H$, $\Vert x\Vert =1$, we
get the desired result (\ref{drag0}).
\end{proof}

\begin{theorem}
Let $T\in B({\mathcal{H}})$ be an $(\alpha ,\beta )$-normal operator, $r\geq
0$ and $\lambda \in {\mathbb{C\diagdown }}\left\{ 0\right\} $. If $\Vert
\lambda T^{\ast }-T\Vert \leq r$ and $\frac{r}{|\lambda |}\leq \inf \{\Vert
T^{\ast }x\Vert :\Vert x\Vert =1\}$, then
\begin{equation}
\alpha ^{2}\Vert T\Vert ^{4}\leq w(T^{2})^{2}+\frac{r^{2}}{|\lambda |^{2}}%
\Vert T\Vert ^{2}.  \label{t4}
\end{equation}
\end{theorem}

\begin{proof}
We use the following reverse of the Schwarz inequality obtained by Dragomir
in \cite{SSD4} (see also \cite[p. 20]{DRABook}):
\begin{equation}
\left( 0\leq \right) \Vert y\Vert ^{2}\,\Vert a\Vert ^{2}-[\text{Re}\langle
y,a\rangle ]^{2}\leq r^{2}\Vert y\Vert ^{2},  \label{r}
\end{equation}%
provided $\Vert y-a\Vert \leq r\leq \Vert a\Vert $.

By the assumption of theorem $\Vert Tx-\lambda T^{\ast }x\Vert \leq r\leq
\Vert \lambda T^{\ast }x\Vert $. Setting $a=\lambda T^{\ast }x$ and $y=Tx$,
with $\Vert x\Vert =1$ in (\ref{r}) we get
\begin{equation*}
\Vert Tx\Vert ^{2}\,\Vert \lambda T^{\ast }x\Vert ^{2}\leq \lbrack \text{Re}%
\langle Tx,\lambda T^{\ast }x\rangle ]^{2}+r^{2}\Vert Tx\Vert ^{2}
\end{equation*}%
whence
\begin{equation}
\alpha ^{2}|\lambda |^{2}\Vert Tx\Vert ^{4}\leq |\lambda |^{2}|\langle
T^{2}x,x\rangle |^{2}+r^{2}\Vert Tx\Vert ^{2}.  \label{rr}
\end{equation}%
Taking the supremum in (\ref{rr}) over $x\in H$, $\Vert x\Vert =1$, we get
the desired result (\ref{t4}).
\end{proof}

Finally, the following result that is less restrictive for the involved
parameters $r$ and $\lambda $ (from the above theorem) may be stated as well:

\begin{theorem}
Let $T\in B({\mathcal{H}})$ be an $(\alpha ,\beta )$-normal operator, $r\geq
0$ and $\lambda \in {\mathbb{C\diagdown }}\left\{ 0\right\} $. If $\Vert
\lambda T^{\ast }-T\Vert \leq r,$ then%
\begin{equation}
\alpha \Vert T\Vert ^{2}\leq w(T^{2})+\frac{r^{2}}{2|\lambda |}.  \label{rra}
\end{equation}
\end{theorem}

\begin{proof}
We use the following reverse of the Schwarz inequality obtained by Dragomir
in \cite{SSD3} (see also \cite[p. 27]{DRABook}):%
\begin{equation}
\left( 0\leq \right) \Vert y\Vert \,\Vert a\Vert -\text{Re}\langle
y,a\rangle \leq \frac{1}{2}r^{2},  \label{rrr}
\end{equation}%
provided $\Vert y-a\Vert \leq r$.

Setting $a=\lambda T^{\ast }x$ and $y=Tx$, with $\Vert x\Vert =1$ in (\ref%
{rrr}) we get
\begin{equation*}
\Vert Tx\Vert \,\Vert \lambda T^{\ast }x\Vert \leq \left\vert \langle
Tx,\lambda T^{\ast }x\rangle \right\vert +\frac{1}{2}r^{2}
\end{equation*}%
which gives%
\begin{equation*}
\alpha \Vert Tx\Vert ^{2}\leq |\langle T^{2}x,x\rangle |+\frac{1}{2|\lambda |%
}r^{2}.
\end{equation*}%
Now, taking the supremum over $\Vert x\Vert =1$ in this inequality, we get
the desired result (\ref{rra})
\end{proof}

\section{Inequalities Involving Norms}

Our first result in this section reads as follows.

\begin{theorem}
Let $T\in B({\mathcal{H}})$ be an $(\alpha ,\beta )$-normal operator. If $%
p\geq 2,$ then
\begin{equation}
2(1+\alpha ^{p})\Vert T\Vert ^{p}\leq \frac{1}{2}(\Vert T+T^{\ast }\Vert
^{p}+\Vert T-T^{\ast }\Vert ^{p}).  \label{2.19b}
\end{equation}

In general, for each $T\in B({\mathcal{H}})$ and $p\geq 2$ we have
\begin{equation}
\left\Vert \frac{T^{\ast }T+TT^{\ast }}{2}\right\Vert ^{p/2}\leq \frac{1}{4}%
(\Vert T+T^{\ast }\Vert ^{p}+\Vert T-T^{\ast }\Vert ^{p}).  \label{2.19bb}
\end{equation}
\end{theorem}

\begin{proof}
We use the following inequality obtained by Dragomir and S\'{a}ndor in \cite%
{DS} (see also \cite[p. 544]{MPF}):
\begin{equation}
\Vert a+b\Vert ^{p}+\Vert a-b\Vert ^{p}\geq 2(\Vert a\Vert ^{p}+\Vert b\Vert
^{p})  \label{2.20b}
\end{equation}%
for any $a,b\in H$ and $p\geq 2$.

Now, if we choose $a=Tx,$ $b=T^{\ast }x$ in (\ref{2.20b}), then we get
\begin{equation}
\Vert Tx+T^{\ast }x\Vert ^{p}+\Vert Tx-T^{\ast }x\Vert ^{p}\geq 2(\Vert
Tx\Vert ^{p}+\Vert T^{\ast }x\Vert ^{p}),  \label{2.21b}
\end{equation}%
whence
\begin{equation}
\Vert Tx+T^{\ast }x\Vert ^{p}+\Vert Tx-T^{\ast }x\Vert ^{p}\geq 2(\Vert
Tx\Vert ^{p}+\alpha ^{p}\Vert Tx\Vert ^{p}),  \label{p2}
\end{equation}%
for any $x\in H,$ $\Vert x\Vert =1.$

Taking the supremum in (\ref{p2}) over $x\in H,$ $\Vert x\Vert =1,$ we get
the desired result (\ref{2.19b}).

Now for the general case $T\in B({\mathcal{H}})$, observe that
\begin{equation}
\Vert Tx\Vert ^{p}+\Vert T^{\ast }x\Vert ^{p}=(\Vert Tx\Vert ^{2})^{\frac{p}{%
2}}+(\Vert T^{\ast }x\Vert ^{2})^{\frac{p}{2}}  \label{2.21bbb}
\end{equation}%
and by applying the elementary inequality:
\begin{equation*}
\frac{a^{q}+b^{q}}{2}\geq \left( \frac{a+b}{2}\right) ^{q},\quad a,b\geq 0%
\text{ \ and \ }q\geq 1
\end{equation*}%
we have
\begin{multline}
(\Vert Tx\Vert ^{2})^{\frac{p}{2}}+(\Vert T^{\ast }x\Vert ^{2})^{\frac{p}{2}%
}\geq 2^{1-\frac{p}{2}}(\Vert Tx\Vert ^{2}+\Vert T^{\ast }x\Vert ^{2})^{%
\frac{p}{2}}  \label{2.11c} \\
=2^{1-\frac{p}{2}}[\langle Tx,Tx\rangle +\langle T^{\ast }x,T^{\ast
}x\rangle ]^{\frac{p}{2}} \\
=2^{1-\frac{p}{2}}[\langle (T^{\ast }T+TT^{\ast })x,x\rangle ]^{\frac{p}{2}}.
\end{multline}%
Combining (\ref{2.21b}) with (\ref{2.11c}) and (\ref{2.21bbb}) we get
\begin{equation}
\frac{1}{4}[\Vert Tx-T^{\ast }x\Vert ^{p}+\Vert Tx+T^{\ast }x\Vert ^{p}]\geq
\left\vert \left\langle \left( \frac{T^{\ast }T+TT^{\ast }}{2}\right)
x,x\right\rangle \right\vert ^{p/2}  \label{2.12c}
\end{equation}%
for any $x\in H,$ $\Vert x\Vert =1$. Taking the supremum over $x\in H,$ $%
\Vert x\Vert =1,$ and taking into account that
\begin{equation*}
w\left( \frac{T^{\ast }T+TT^{\ast }}{2}\right) =\left\Vert \frac{T^{\ast
}T+TT^{\ast }}{2}\right\Vert ,
\end{equation*}%
we deduce the desired result (\ref{2.19bb}).
\end{proof}

\begin{theorem}
Let $T\in B({\mathcal{H}})$ be an $(\alpha ,\beta )$-normal operator. If $%
p\in (1,2)$ and $\lambda ,\mu \in {\mathbb{C}}$, then
\begin{multline}
\left[ (|\lambda |+\beta |\mu |)^{p}+\max \{|\lambda |-|\mu |\beta ,\alpha
|\mu |-|\lambda |\}\right] \Vert T\Vert ^{p}  \label{lm} \\
\leq \Vert \lambda T+\mu T^{\ast }\Vert ^{p}+\Vert \lambda T-\mu T^{\ast
}\Vert ^{p}.
\end{multline}
\end{theorem}

\begin{proof}
We use the following inequality obtained by Dragomir and S\'{a}ndor in \cite%
{DS} (see also \cite[p. 544]{MPF})
\begin{equation}
(\Vert a\Vert +\Vert b\Vert )^{p}+\left\vert \,\Vert a\Vert -\Vert b\Vert
\,\right\vert ^{p}\leq ||a+b\Vert ^{p}+\Vert a-b\Vert ^{p},  \label{p}
\end{equation}%
for any $a,b\in H$ and $p\in (1,2)$.

Put $a=\lambda Tx,$ $b=\mu T^{\ast }x$ in (\ref{p}) to obtain
\begin{multline*}
(\Vert \lambda Tx\Vert +\Vert \mu T^{\ast }x\Vert )^{p}+\left\vert \Vert
\lambda Tx\Vert -\Vert \mu T^{\ast }x\Vert \right\vert ^{p} \\
\leq \Vert \lambda Tx+\mu T^{\ast }x\Vert ^{p}+\Vert \lambda Tx-\mu T^{\ast
}x\Vert ^{p},
\end{multline*}%
whence
\begin{multline}
(|\lambda |+|\mu |\alpha )^{p}\Vert Tx\Vert ^{p}+\left( \max \{|\lambda
|-|\mu |\beta ,\alpha |\mu |-|\lambda |\}\right) \Vert Tx\Vert ^{p}
\label{pp} \\
\leq \Vert \lambda Tx+\mu T^{\ast }x\Vert ^{p}+\Vert \lambda Tx-\mu T^{\ast
}x\Vert ^{p},
\end{multline}%
for any $x\in H,$ $\Vert x\Vert =1.$

Taking the supremum in (\ref{pp}) over $x\in H,$ $\Vert x\Vert =1,$ we get
the desired result (\ref{lm}).
\end{proof}

\textbf{Acknowledgement.} The authors would like to thank Professor
Leiba Rodman for his useful suggestions and for bringing the
majorization lemma of R.G. Douglas into their attention. The authors
would also like to express their gratitude to Professor M.
Mirzavaziri for giving an example of an $(\alpha ,\beta )$-normal
operator which is neither normal nor hyponormal.

\end{document}